\pdfoutput=1
\RequirePackage{ifpdf}
\ifpdf % We~are running pdfTeX in pdf mode
\documentclass[pdftex]{sigma}
\else
\documentclass{sigma}
\fi

\numberwithin{equation}{section}

\newtheorem*{Theorem*}{Theorem}

 { \theoremstyle{definition}

 }
\newcommand{\qbinomial}[2]{\genfrac{[}{]}{0pt}{}{#1}{#2}}

\begin{document}

%\allowdisplaybreaks

\renewcommand{\thefootnote}{}

\newcommand{\arXivNumber}{2408.00377}

\renewcommand{\PaperNumber}{103}

\FirstPageHeading

\ShortArticleName{Rogers--Ramanujan Type Identities Involving Double Sums}

\ArticleName{Rogers--Ramanujan Type Identities\\ Involving Double Sums\footnote{This paper is a~contribution to the Special Issue on Basic Hypergeometric Series Associated with Root Systems and Applications in honor of Stephen C.~Milne's 75th birthday. The~full collection is available at \href{https://www.emis.de/journals/SIGMA/Milne.html}{https://www.emis.de/journals/SIGMA/Milne.html}}}

\Author{Dandan CHEN~$^{\rm ab}$ and Siyu YIN~$^{\rm a}$}

\AuthorNameForHeading{D.~Chen and S.~Yin}

\Address{$^{\rm a)}$~Department of Mathematics, Shanghai University, P.R. China}
\EmailD{\href{mailto:mathcdd@shu.edu.cn}{mathcdd@shu.edu.cn}, \href{mailto:siyuyin@shu.edu.cn}{siyuyin@shu.edu.cn}}
\URLaddressD{\url{https://math.shu.edu.cn/info/1021/3135.htm}}

\Address{$^{\rm b)}$~Newtouch Center for Mathematics, Shanghai University, P.R. China}

\ArticleDates{Received August 02, 2024, in final form November 13, 2024; Published online November 19, 2024}

\Abstract{We prove four new Rogers--Ramanujan-type identities for double series. They follow from the classical Rogers--Ramanujan identities using the constant term method and properties of Rogers--Szeg\H{o} polynomials.}


\Classification{11P84; 33D15}

\renewcommand{\thefootnote}{\arabic{footnote}}
\setcounter{footnote}{0}

\section{Introduction}
In 1894, L.J.~Rogers \cite{Rogers-1894} discovered numerous sum-product $q$-series identities. Among his findings, he proved the following identities:
\begin{align}
&\sum_{n=0}^{\infty}\frac{q^{n^2}}{(q;q)_n}=\frac{1}{\bigl(q,q^4;q^5\bigr)_{\infty}},\label{1-4-5}\\
&\sum_{n=0}^{\infty}\frac{q^{n^2+n}}{(q;q)_n}=\frac{1}{\bigl(q^2,q^3;q^5\bigr)_{\infty}},\label{2-3-5}\\
&\sum_{n=0}^{\infty}\frac{q^{n^2}}{\bigl(q^4;q^4\bigr)_n}=\frac{1}{\bigl(-q^2;q^2\bigr)_{\infty}\bigl(q,q^4;q^5\bigr)_{\infty}},\label{n1-4-5}\\
&\sum_{n=0}^{\infty}\frac{q^{n^2+2n}}{\bigl(q^4;q^4\bigr)_n}=\frac{1}{\bigl(-q^2;q^2\bigr)_{\infty}\bigl(q^2,q^3;q^5\bigr)_{\infty}}.\label{n2-3-5}
\end{align}

We briefly introduce the notations used in this paper. We always assume $|q|< 1$ for convergence. The standard $q$-series notations are as follows \cite{Gasper-Rahman-2004}:
\begin{align*}
&(a;q)_0:=1,\qquad (a;q)_n:=\prod_{k=0}^{n-1}\bigl(1-aq^k\bigr),\qquad (a;q)_{\infty}:=\prod_{k=0}^{\infty}\bigl(1-aq^k\bigr),\\
&(a_1,\dots,a_m;q)_n=(a_1;q)_n\cdots(a_m;q)_n,\qquad n\in\mathbb{N}\cup\{\infty\}.
\end{align*}

The identities referred to as \eqref{1-4-5} and \eqref{2-3-5} are recognized as the Rogers--Ramanujan identities, having been rediscovered by Ramanujan prior to the year 1913. They have inspired a lot of work on finding identities of similar forms. Rogers--Ramanujan type identities function as one of the witnesses for deep connections between the theory of $q$-series and modular forms. After multiplying with suitable powers of $q$, the right-hand side of \eqref{1-4-5} and \eqref{2-3-5} become modular forms which is not clear from the sum sides. An important question in the theory of $q$-series and modular forms is to judge what kind of basic hypergeometric series qualify as modular forms. This question remains an open challenge in the field. In a series of works, W.~Nahm \cite{Nahm-1994, Nahm-1995, Nahm-2007} considered the series
\begin{align*}
f_{A,B,C}(q):=\sum_{n=(n_1,\dots,n_r)^T\in(\mathbb{Z}_{\geq 0})^r}\frac{q^{\frac{1}{2}n^TAn+n^TB+C}}{(q;q)_{n_1}\cdots(q;q)_{n_r}},
\end{align*}
where $r\geq 1$ is a positive integer, $A$ is a real positive definite symmetric $r\times r$ matrix, $B$ is a~vector of length $r$, and $C$ is a scalar.

Nahm \cite{Nahm-2007} proposed a conjecture that provides sufficient and necessary conditions on the matrix part of a modular triple. The conjecture is formulated in terms of the Bloch group and a system of polynomial equations induced by the matrix part. D.~Zagier \cite{Zagier-2007} gave a precise statement of this conjecture. When the rank $r=1$, the identities \eqref{1-4-5}--\eqref{n2-3-5} showed that\looseness=1
\begin{align*}
(A,B,C)=(2,0,-1/60),\qquad(2,1,11/60), \qquad(1/2,0,-1/40), \qquad(1/2,1/2,1/40)
\end{align*}
are all modular triples. Zagier \cite{Zagier-2007} studied Nahm's problem and identified many possible modular triples. In particular, for rank $r=1$, Zagier substantiated Nahm's conjecture and proved that there exactly seven modular triples. Besides the four aforementioned triples, the remaining triples are
\begin{align*}
(1,0,-1/48),\qquad(1,1/2,1/24),\qquad(1,-1/2,1/24),
\end{align*}
which is easily justified by Euler's identities.

Following the notion in \cite{Wang-2023}, some Rogers--Ramanujan type identities are characterized by a~distinct structural pattern. For a given integer $k$, an identity of the following shape is defined as finite sum of
\[
\sum_{(i_1,\dots,i_k)\in S}\frac{(-1)^{t(i_1,\dots,i_k)}q^{Q(i_1,\dots,i_k)}}{(q^{n_1};q^{n_1})_{i_1}\cdots(q^{n_k};q^{n_k})_{i_k}}=\prod_{(a,n)\in P}(q^a;q^n)_{\infty}^{r(a,n)}
\]
as  Rogers--Ramanujan type identities of $\operatorname{index}(n_1,n_2,\dots,n_k)$. Here $t(i_1,\dots,i_k)$ is an integer-valued function, $Q(i_1,\dots,i_k)$ is a rational polynomials in variables $i_1,\dots,i_k,n_1,\dots,n_k$ are positive integers with $\gcd(n_1,n_2,\dots,n_k)=1$, $S$ is a subset of $\mathbb{Z}^k$, $P$ is a finite subset of $\mathbb{Q}^2$ and $r(a,n)$ are integer-valued functions.

In 2021, Andrews and Uncu \cite{Andrews-Uncu-2023} proved an identity of index $(1,3)$ and further conjectured that \cite[Conjecture 1.2]{Andrews-Uncu-2023}
\begin{align*}
\sum_{i,j\geq 0}\frac{(-1)^jq^{3j(3j+1)/2+i^2+3ij+i+j}}{(q;q)_i\bigl(q^3;q^3\bigr)_j}=\frac{1}{\bigl(q^2,q^3;q^6\bigr)_{\infty}}.
\end{align*}
This was first proved by Chern \cite{Chern-2022} and then by Wang \cite{Wang-2023}. Besides, Cao and Wang \cite{MR4530918} established some Rogers--Ramanujan type identities of indexes
\[
(1,1),\qquad(1,2),\qquad(1,1,1),\qquad(1,1,3),\qquad(1,2,2),\qquad(1,2,3),\qquad(1,2,4).
\]
For instance, they proved that for any $u\in \mathbb{C}$
\begin{align*}
&\sum_{i,j,k\geq 0}\frac{(-1)^{i+j}u^{i+3k}q^{(i^2-i)/2+(i-2j+3k)^2/4}}{(q;q)_i\bigl(q^2;q^2\bigr)_j\bigl(q^3;q^3\bigr)_k} =\frac{\bigl(u^2;q\bigr)_{\infty}\bigl(q,-u^2;q^2\bigr)_{\infty}}{\bigl(-u^6;q^6\bigr)_{\infty}},\\
&\sum_{i,j,k\geq 0}\frac{(-1)^{(i-2j+3k)/2}u^{i+k}q^{(i^2-i)/2+(i-2j+3k)^2/4}}{(q;q)_i\bigl(q^2;q^2\bigr)_j\bigl(q^3;q^3\bigr)_k} =\frac{\bigl(q;q^2\bigr)_{\infty}\bigl(-u^2;q^3\bigr)_{\infty}}{\bigl(u^2;q^6\bigr)_{\infty}}.
\end{align*}
Furthermore, Wang proved Zagier's rank three examples for Nahm's problem one by one in \cite{Wang-2024}.

Motivated by the constant term method \cite{Andrews-1986} and the identities \eqref{1-4-5}--\eqref{n2-3-5}, we present the following theorem.
\begin{theorem}\label{main-thm}
We have
\begin{align}
&\sum_{n,m\geq 0}\frac{(-1)^{\binom{n-m}{2}}q^{\frac{3m^2}{4}+\frac{mn}{2}+\frac{3n^2}{4}}}{(q;q)_m(q;q)_n}=\frac{1}{\bigl(q^2,q^8;q^{10}\bigr)_{\infty}},\label{2-8-10}\\
&\sum_{n,m\geq 0}\frac{(-1)^{\binom{n-m}{2}}q^{\frac{3m^2}{4}+\frac{mn}{2}+\frac{3n^2}{4}+m+n}}{(q;q)_m(q;q)_n}=\frac{1}{\bigl(q^4,q^6;q^{10}\bigr)_{\infty}},\label{4-6-10}\\
&\sum_{n,m\geq 0}\frac{(-1)^mq^{\frac{m^2}{4}+\frac{mn}{2}+\frac{n^2}{4}}}{\bigl(q^2;q^2\bigr)_m\bigl(q^2;q^2\bigr)_n} =\frac{1}{\bigl(-q^2;q^2\bigr)_{\infty}\bigl(q,q^4;q^5\bigr)_{\infty}},\label{n2-1-4-5}\\
&\sum_{n,m\geq 0}\frac{(-1)^mq^{\frac{m^2}{4}+\frac{mn}{2}+\frac{n^2}{4}+m+n}}{\bigl(q^2;q^2\bigr)_m\bigl(q^2;q^2\bigr)_n} =\frac{1}{\bigl(-q^2;q^2\bigr)_{\infty}\bigl(q^2,q^3;q^5\bigr)_{\infty}}\label{n2-2-3-5}.
\end{align}
\end{theorem}

The paper is organized as follows. In Section \ref{sec-pre}, we introduce some basic identities and the constant term method. In Section \ref{proof}, we demonstrate the proof of Theorem \ref{main-thm} using the constant term method.

\section{Preliminaries}\label{sec-pre}
In this section, we first collect some useful identities on basic hypergeometric series. The $q$-bino\-mi\-al theorem \cite[p.~8]{Gasper-Rahman-2004} is defined as
\begin{align*}
\frac{(az;q)_\infty}{(z;q)_\infty}
=\sum_{n=0}^{\infty}\frac{(a;q)_n}{(q;q)_n}z^n,\qquad |z|<1.
\end{align*}

As corollaries, Euler's $q$-exponential identities assert \cite{Gasper-Rahman-2004}:
\begin{align*}
&\sum_{n=0}^{\infty}\frac{z^n}{(q;q)_n}=\frac{1}{(z;q)_{\infty}},\qquad|z|<1,\\
&\sum_{n=0}^{\infty}\frac{q^{\binom{n}{2}}z^n}{(q;q)_n}=(-z;q)_{\infty}.
\end{align*}

The Jacobi's triple product identity is given by \cite[p.~15]{Gasper-Rahman-2004}
\begin{align*}
\bigl(q,zq^{\frac12},q^{\frac12}/z;q\bigr)_{\infty}=\sum_{n=-\infty}^{\infty}(-1)^nq^{n^2/2}z^n,\qquad z\neq0.
\end{align*}

Recall the new representation of Rogers--Szeg\H{o} polynomials given by Berkovich and Warnaar~\mbox{\cite[Theorem 8.1]{Berkovich-Warnaar-2005}}
\begin{align}\label{RS-polynomial}
H_n(t;q)=\sum_{r=0}^{\lfloor \frac{n}{2}\rfloor}t^{2r}\bigl(-q/t;q^2\bigr)_r\bigl(-t;q^2\bigr)_{\lfloor \frac{n+1}{2}\rfloor-r}
\qbinomial{\lfloor n/2\rfloor}{r}_{q^2},
\end{align}
where $H_n(t;q)$ was originally defined as
\begin{align*}
H_n(t;q)=\sum_{j=0}^n t^j\qbinomial{n}{j}_{q}.
\end{align*}
By \eqref{RS-polynomial}, we directly deduce that $H_{2n}(-1;q)=\bigl(q;q^2\bigr)_n$ and $H_{2n+1}(-1;q)=0$, which are also deduced from the generating function for these polynomials, that was known to Rogers \cite{Gasper-Rahman-2004}.

Besides, for a series $f(z)=\sum_{n=-\infty}^{\infty}a(n)z^n$, we denote the coefficient of $z^n$ by $[z^n]f(z)=a(n)$. Specifically, we use $CT_zf(z)$ to denote the constant term $\bigl[z^0\bigr]f(z)$. It is a well-established fact~that
\begin{align*}
\oint_K f(z)\frac{{\rm d}z}{2\pi {\rm i}z}=CT_zf(z)=\bigl[z^0\bigr]f(z),
\end{align*}
where $K$ is a positively oriented, simple closed contour around the origin. Consequently, we often calculate the constant term or integral. By doing so we can transform the original series into new series, which can then be assessed by some well-known identities. For simplicity, when integral calculus is not necessary, we prefer to utilize the constant term method in proofs.

\section{The proof of Theorem~\ref{main-thm}}\label{proof}
In this section, we present the proof of Theorem \ref{main-thm} by the constant term method and famous $q$-series identities.
\begin{proof}
For \eqref{2-8-10}, we obtain that
\begin{align*}
&{}\sum_{n,m\geq 0}\frac{(-1)^{\binom{n-m}{2}}q^{\frac{3}{4}m^2+\frac{1}{2}mn+\frac{3}{4}n^2}}{(q;q)_m(q;q)_n}
=\sum_{n,m\geq 0}\frac{{\rm i}^{n-m}q^{\frac{3}{4}m^2+\frac{1}{2}mn+\frac{3}{4}n^2}}{(q;q)_m(q;q)_n}\\
&\qquad{}=\sum_{n,m\geq 0}\frac{{\rm i}^{n-m}q^{\frac{1}{2}\binom{m+n}{2}+\binom{m}{2}+\frac{3m}{4}+\binom{n}{2}+\frac{3n}{4}}}{(q;q)_m(q;q)_n}\\
&\qquad{}=CT_{z}\sum_{m\geq 0}\frac{{\rm i}^mz^mq^{\binom{m}{2}+\frac{3m}{4}}}{(q;q)_m}\sum_{n\geq 0}\frac{(-{\rm i})^nz^nq^{\binom{n}{2}+\frac{3n}{4}}}{(q;q)_n}
\sum_{k=-\infty}^{\infty}(-1)^kq^{\frac{1}{2}\binom{k}{2}}z^{-k}\\
&\qquad{}=CT_{z}\bigl(-z^2q^{\frac{3}{2}};q^2\bigr)_{\infty}\sum_{k=-\infty}^{\infty}(-1)^kq^{\frac{1}{2}\binom{k}{2}}z^{-k}\\
&\qquad{}=CT_{z}\sum_{n\geq 0}\frac{q^{2\binom{n}{2}}z^{2n}q^{\frac{3n}{2}}}{\bigl(q^2;q^2\bigr)_n}\sum_{k=-\infty}^{\infty}(-1)^kq^{\frac{1}{2}\binom{k}{2}}z^{-k}\\
&\qquad{}=\sum_{n=0}^{\infty}\frac{q^{2n^2}}{\bigl(q^2;q^2\bigr)_n}.
\end{align*}
Using \eqref{1-4-5}, we complete the proof of \eqref{2-8-10}.

For \eqref{4-6-10}, we similarly write the item \smash{$q^{\frac{1}{2}\binom{m+n}{2}}$} into \smash{$q^{\frac{1}{2}\binom{k}{2}}$}. By the constant item method and Euler's $q$-exponential identities, we can easily prove it.

After exchanging the summation, we prove \eqref{n2-1-4-5} using the special case of \eqref{RS-polynomial},{\samepage
\begin{align*}
\sum_{n,m\geq 0}\frac{(-1)^mq^{\frac{m^2}{4}+\frac{mn}{2}+\frac{n^2}{4}}}{\bigl(q^2;q^2\bigr)_m\bigl(q^2;q^2\bigr)_n}
&=\sum_{n\geq 0}\sum_{m=0}^n \frac{(-1)^nq^{\frac{n^2}{4}}}{\bigl(q^2;q^2\bigr)_n}\cdot(-1)^m\qbinomial{n}{m}_{q^2}\\
&=\sum_{n\geq 0} \frac{(-1)^nq^{\frac{n^2}{4}}}{\bigl(q^2;q^2\bigr)_n}\cdot H_n\bigl(-1;q^2\bigr)\\
&=\sum_{n\geq 0}\frac{(-1)^{2n}q^{n^2}}{\bigl(q^2;q^2\bigr)_{2n}}\cdot \bigl(q^2;q^4\bigr)_{n}
=\sum_{n\geq 0}\frac{q^{n^2}}{\bigl(q^4;q^4\bigr)_n}.
%&=\frac{1}{(-q^2;q^2)_{\infty}(q,q^4;q^5)_{\infty}}.
\end{align*}
With the help of \eqref{n1-4-5} and \eqref{n2-1-4-5} is proved.}

Similar to the method used in \eqref{n2-1-4-5} and \eqref{n2-2-3-5} can also be proved by the constant term method,
\begin{align*}
\sum_{n,m\geq 0}\frac{(-1)^mq^{\frac{m^2}{4}+\frac{mn}{2}+\frac{n^2}{4}+m+n}}{\bigl(q^2;q^2\bigr)_m\bigl(q^2;q^2\bigr)_n}
&=\sum_{n,m\geq 0}\frac{(-{\rm i})^{n-m}q^{\frac{3m+3n}{2}}}{\bigl(q^2;q^2\bigr)_m\bigl(q^2;q^2\bigr)_n}\cdot q^{\frac{(m+n)(m+n-2)}{4}}{\rm i}^{n+m}\\
&=CT_z\sum_{m\geq 0}\frac{q^{\frac{3m}{2}}z^{m}(i)^m}{\bigl(q^2;q^2\bigr)_m}\sum_{n\geq 0}\frac{q^{\frac{3n}{2}}z^{n}(-{\rm i})^n} {\bigl(q^2;q^2\bigr)_n}\sum_{k=-\infty}^{\infty}q^{\frac{k(k-2)}{4}}z^{-k}{\rm i}^k\\
&=CT_z\frac{1}{\bigl(-z^2q^3;q^4\bigr)_{\infty}}\sum_{k=-\infty}^{\infty}q^{\frac{k(k-2)}{4}}z^{-k}{\rm i}^k\\
&=\sum_{n=0}^{\infty}\frac{\bigl(-z^2q^3\bigr)^n}{\bigl(q^4;q^4\bigr)_n}(-1)^nz^{-2n}q^{2\binom{n}{2}}\\
&=\sum_{n=0}^{\infty}\frac{q^{n^2+2n}}{\bigl(q^4;q^4\bigr)_n}.
%&=\frac{1}{(-q^2;q^2)_{\infty}(q^2,q^3;q^5)_{\infty}}.
\end{align*}
Combining with \eqref{n2-3-5}, we find that \eqref{n2-2-3-5} holds.
\end{proof}

\subsection*{Acknowledgements}
We are grateful to the referees and the editors for their helpful comments and suggestions. We thank Warnaar for some valuable comments, especially for bringing the work \cite{Berkovich-Warnaar-2005} to our attention.
The first author was supported in part by the National Natural Science Foundation of China (Grant No. 12201387).

%\bibliographystyle{sigma}
%\bibliography{example}

\begin{thebibliography}{99}
\footnotesize\itemsep=0pt

\bibitem{Andrews-1986}
Andrews G.E., {$q$}-series: their development and application in analysis,
 number theory, combinatorics, physics, and computer algebra, \textit{CBMS
 Reg. Conf. Ser. Math.}, Vol.~66, \href{https://doi.org/10.1090/cbms/066}{American Mathematical Society}, Providence,
 RI, 1986.

\bibitem{Andrews-Uncu-2023}
Andrews G.E., Uncu A.K., Sequences in overpartitions, \href{https://doi.org/10.1007/s11139-022-00685-y}{\textit{Ramanujan~J.}}
 \textbf{61} (2023), 715--729, \href{https://arxiv.org/abs/2111.15003}{arXiv:2111.15003}.

\bibitem{Berkovich-Warnaar-2005}
Berkovich A., Warnaar S.O., Positivity preserving transformations for
 {$q$}-binomial coefficients, \href{https://doi.org/10.1090/S0002-9947-04-03680-3}{\textit{Trans. Amer. Math. Soc.}} \textbf{357}
 (2005), 2291--2351, \href{https://arxiv.org/abs/math.CO/0302320}{arXiv:math.CO/0302320}.

\bibitem{MR4530918}
Cao Z., Wang L.,
Multi-sum Rogers--Ramanujan type identities,
\textit{J.~Math. Anal. Appl.} \textbf{522} (2023), 126960, 24~pages, \href{https://arxiv.org/abs/2205.12786}{arXiv:2205.12786}.

\bibitem{Chern-2022}
Chern S., Asymmetric {R}ogers--{R}amanujan type identities.~{I}. {T}he
 {A}ndrews--{U}ncu conjecture, \href{https://doi.org/10.1090/proc/16332}{\textit{Proc. Amer. Math. Soc.}} \textbf{151}
 (2023), 3269--3279, \href{https://arxiv.org/abs/2203.15168}{arXiv:2203.15168}.

\bibitem{Gasper-Rahman-2004}
Gasper G., Rahman M., Basic hypergeometric series, 2nd ed., \textit{Encyclopedia Math.
 Appl.}, Vol.~96, \href{https://doi.org/10.1017/CBO9780511526251}{Cambridge University Press}, Cambridge, 2004.

\bibitem{Nahm-1994}
Nahm W., Conformal field theory, dilogarithms, and three-dimensional manifolds,
 \textit{Adv. Appl. Clifford Algebras} \textbf{4} (1994), 179--191.

\bibitem{Nahm-1995}
Nahm W., Conformal field theory and the dilogarithm, in X{I}th {I}nternational
 {C}ongress of {M}athematical {P}hysics ({P}aris, 1994), International Press,
 Cambridge, MA, 1995, 662--667.

\bibitem{Nahm-2007}
Nahm W., Conformal field theory and torsion elements of the {B}loch group, in
 Frontiers in Number Theory, Physics, and Geometry.~{II}, \href{https://doi.org/10.1007/978-3-540-30308-4_2}{Springer}, Berlin,
 2007, 67--132, \href{https://arxiv.org/abs/hep-th/0404120}{arXiv:hep-th/0404120}.

\bibitem{Rogers-1894}
Rogers L.J., Second memoir on the expansion of certain infinite products,
 \href{https://doi.org/10.1112/plms/s1-25.1.318}{\textit{Proc. Lond. Math. Soc.}} \textbf{25} (1893), 318--343.

\bibitem{Wang-2023}
Wang L., New proofs of some double sum {R}ogers--{R}amanujan type identities,
 \href{https://doi.org/10.1007/s11139-022-00654-5}{\textit{Ramanujan~J.}} \textbf{62} (2023), 251--272, \href{https://arxiv.org/abs/2203.15572}{arXiv:2203.15572}.

\bibitem{Wang-2024}
Wang L., Explicit forms and proofs of {Z}agier's rank three examples for
 {N}ahm's problem, \href{https://doi.org/10.1016/j.aim.2024.109743}{\textit{Adv. Math.}} \textbf{450} (2024), 109743, 46~pages,
 \href{https://arxiv.org/abs/2211.04375}{arXiv:2211.04375}.

\bibitem{Zagier-2007}
Zagier D., The dilogarithm function, in Frontiers in Number Theory, Physics,
 and Geometry.~{II}, \href{https://doi.org/10.1007/978-3-540-30308-4_1}{Springer}, Berlin, 2007, 3--65.

\end{thebibliography}

\pdfbookmark[1]{References}{ref}
\LastPageEnding

\end{document}